\documentclass[leqno, 12pt]{article}

\usepackage{amsmath}
\usepackage{amsfonts}
\usepackage{amssymb}
\usepackage{amsthm} 
\usepackage{graphicx}

\usepackage{scalerel,stackengine}
\stackMath
\newcommand\widecheck[1]{%
\savestack{\tmpbox}{\stretchto{%
  \scaleto{%
    \scalerel*[\widthof{\ensuremath{#1}}]{\kern-.6pt\bigwedge\kern-.6pt}%
    {\rule[-\textheight/2]{1ex}{\textheight}}
  }{\textheight}%
}{0.5ex}}%
\stackon[1pt]{#1}{\scalebox{-1}{\tmpbox}}%
}
\usepackage{scalerel,stackengine}
\stackMath
\newcommand\widebreve[1]{%
\savestack{\tmpbox}{\stretchto{%
  \scaleto{%
    \scalerel*[\widthof{\ensuremath{#1}}]{\kern-.1pt\bigcup\kern-.1pt}  
    {\rule[-\textheight/2]{1ex}{\textheight}}
  }{\textheight}%
}{0.5ex}}%
\stackon[1pt]{#1}{\scalebox{1}{\tmpbox}}%
}
\usepackage{comment}

\def\tod{\buildrel d\over\to}
\def\toL1{\buildrel \mathit{L_1}\over\longrightarrow}

\def\towpr1{\buildrel w.pr.1\over\to}
\def\11{\buildrel 1-1\over\longleftrightarrow}
\def\eqd{\buildrel d\over=}

\def\~{\tilde}
\def\stms{\!\setminus\!} 
\def\^{\hat}

\def\c{\check}

\def\implies{\Rightarrow}

\def\Lam{\Lambda}
\def\Gam{\Gamma}
\def\gam{\gamma}
\def\lam{\lambda}
\def\alp{\alpha}
\def\bet{\beta}
\def\del{\delta}
\def\Del{\Delta}
\def\eps{\epsilon}

\def\prtl{\partial}

\def\E{{\rm E}}
\def\P{{\rm P}}

\def\R1{{\bf R}^1}
\def\B1{{\cal B}^1}

\def\nid{\noindent}

\def\ts{\textstyle}
\def\t12{\textstyle{1\over2}}

\def\smalltype{\let\rm=\eightrm \let\bf=\eightbf \let\it=\eightit \let\sl=\eightsl
 \baselineskip=9pt \rm}
\font\tenrm=cmr10
\font\tenbf=cmbx10
\font\tenit=cmti10
\font\tensl=cmsl10
\def\medtype{\let\rm=\tenrm \let\bf=\tenbf \let\it=\tenit \let\sl=\tensl
 \baselineskip=12pt \rm}
\font\twelverm=cmr12
\font\twelvebf=cmbx12
\font\twelveit=cmti12 
\font\twelvesl=cmsl12
\def\bigtype{\let\rm=\twelverm \let\bf=\twelvebf \let\it=\twelveit \let\sl=\twelvesl
 \baselineskip=14pt \rm}
 
\begin{document}

\title{Pure Significance Tests for Multinomial and Binomial Distributions: the Uniform Alternative
{}}
\author{Michael D. Perlman\footnote{mdperlma@uw.edu.}\\Department of Statistics\\
University of Washington\\Seattle WA 98195, U.S.A.}

\maketitle

\begin{abstract}

\nid  A {\it pure significance test} (PST) tests a simple null hypothesis $H_f:Y\sim f$ {\it without specifying an alternative hypothesis} by rejecting $H_f$ for {\it small} values of $f(Y)$.  When the sample space  supports a proper uniform pmf $f_\mathrm{unif}$, the PST can be viewed as a classical likelihood ratio test for testing $H_f$ against this uniform alternative.  Under this interpretation, standard test features such as power, Kullback-Leibler divergence, and expected $p$-value can be considered. This report focuses on PSTs for multinomial and binomial distributions, and for the related goodness-of-fit testing problems with the uniform alternative. The case of repeated observations cannot be reduced to the single observation case via sufficiency. The {\it ordered binomial distribution}, apparently new, arises in the course of this study.

\end{abstract}

\newpage

\nid{\bf 1. Pure Significance Tests.}

\nid Let $Y$ denote a random vector (rvtr) with sample space ${\cal Y}$ and let $f$ be a probability mass function (pmf) $f$ on ${\cal Y}$. A {\it pure significance test (PST)} for testing the simple null hypothesis $H_f:Y\sim f$  {\it without specifying an alternative hypothesis} rejects $H_f$ for small values of $f(Y)$. If $Y=y_0$ is observed, the attained $p$ value is
\begin{equation}\label{PST1}
\P_f[f(Y)\le f(y_0)].
\end{equation}

Pure significance tests long have been a contentious subject, since many believe that relative likelihood, not likelihood, is appropriate for non-Bayesian statistical analysis, hence an alternative hypothesis must be specified. See Hodges (1990) and Howard (2009) for informative reviews and see Appendix 2 of this report for a critique of PSTs. Here we attempt to reconcile, or at least reduce, this contention by noting that when the sample space ${\cal Y}$ supports a proper uniform pmf $f_\mathrm{unif}$, the PST can be viewed as a classical likelihood ratio test (LRT) for testing 
\begin{equation}\label{PST2}
H_f:Y\sim f\quad\mathrm{vs.}\quad H_\mathrm{unif}:Y\sim f_\mathrm{unif},
\end{equation}
the uniform alternative. This interpretation allows consideration of standard test features for PST's, such as power, Kullback-Leibler divergence (KLD), and expected $p$-value (EPV).

This report studies PSTs for multinomial distributions, with binomial distributions as a special case.\footnote{\label{NSASAG}This topic arose in Problem  2023-03 under the UW Statistics Department's NSASAG study project.} After briefly reviewing the geometry of the multinomial family in Section 2, problem \eqref{PST2} for testing a simple multinomial hypothesis against uniformity is discussed in Section 3. It is shown in Propostion 3.1 that under the KLD criterion, the equal-cell-probability (ecp) multinomial is closest to the uniform alternative, which is not a member of  the multinomial family, so that the $p$-value for this case may provide an upper bound for the $p$-value in the general case. The $p$-value is then examined directly by studying the EPV criterion, where Conjecture 3.3 with supporting results provide further evidence for this upper bound property. 

A related testing problem in treated in Section 4, where the simple null hypothesis in \eqref{PST2} is replaced by the composite hypothesis consisting of the entire multinomial family; again the alternative hypothesis is the uniform distribution.  This can be viewed as testing goodness-of-fit for the multinomial family, but with a specified alternative, namely the uniform distribution. Here the ecp multinomial distribution is not the least favorable distribution for this problem (cf. Proposition 4.1), so the PST with the ecp distribution as the null hypothesis is inappropriate for the composite multinomial null hypothesis. Instead, the LRT for the composite hypothesis is derived and shown to be unbiased, with its $p$-value determined by the ecp distribution (cf. Proposition 4.2 and Remark 4.3).

The first testing problem is revisited in Sections 5-7 for the apparently simpler but still challenging binomial case, where Conjecture 3.3 for the EPV criterion is refined and verified in several cases. Hopefully this will lead to further insight for the general multinomial case. The binomial case leads to the introduction of the {\it ordered binomial distribution (OBD)} (Definition 7.1), which appears to be new and of interest in its own right.

Although the preceding results are stated for a single observation from a multinomial or binomial distribution, at first glance it might appear that they extend directly to the case of repeated observations, since the sum of the observations is a sufficient statistic whose distribution remains in the multinomial or binomial family. In Section 8, however, we note that this is invalid for the PST when the uniform alternative is introduced, since the sum statistic is not sufficient for the multinomial/binomial $+$ uniform model. The repeated observations case appears to be substantially more challenging.

Most proofs appear in Appendix 1, while Appendix 2 presents a general critique of pure significance tests. Warm thanks go to Steven W. Knox for helpful discussions.
\vskip6pt

\nid{\bf 2. Some geometry of the multinomial family.}  

\nid  Let $\mathrm{Mult}_{\bf p}(k,n)$ denote the $k$-cell multinomial distribution based on $n$ trials with cell probabilities $(p_1,\dots,p_k)\equiv{\bf p}$. Its sample space  is the integer simplex
\begin{equation}\label{samplespa}
\ts {\cal S}_{k,n}=\{{\bf x}\equiv(x_1,\dots,x_k)\,|\,x_1\ge0,\dots,x_k\ge0,\,\sum_{j=1}^k x_j=n\},
\end{equation}
while its pmf  is given by
\begin{equation}\label{pmfa}
 f_{\bf p}({\bf x})
={n\choose{\bf x}}\prod_{j=1}^k p_j^{x_j},\quad {\bf x}\in{\cal S}_{k,n},
\end{equation}
where ${n\choose{\bf x}}=\frac{n!}{\prod_{j=1}^k x_j!}$ is the multinomial coefficient. Denote the set of all such multinomial pmfs by
\begin{equation}\label{setMult}
{\cal M}(k,n)=\{f_{\bf p}\,|\,{\bf p}\in{\cal P}_k\},
\end{equation}
where ${\cal P}_k$ is the $(k-1)$-dimensional probability simplex:
\begin{equation}\label{psimp}
{\cal P}_k=\{{\bf p}\equiv(p_1,\dots,p_k)\,|\,0\le p_j\le 1,\sum_{j=1}^k p_j=1\}.
\end{equation}

Because $|{\cal S}_{k,n}|={n+k-1\choose k-1}$, the multinomial pmf $f_{\bf p}$ can be viewed as a vector
 in the probability simplex ${\cal P}_{n+k-1\choose k-1}$ and the multinomial family ${\cal M}(k,n)$
is a $(k-1)$-dimensional curved surface in ${\cal P}_{n+k-1\choose k-1}$. For example, in the simplest case $k=n=2$ (binomial, 2 trials), ${\bf p}=(p,1-p)$, $0\le p\le1$, and 
\begin{equation*}
f_{\bf p}=\big((1-p)^2,\,2p(1-p),\,p^2\big),
\end{equation*}
so ${\cal M}(2,2)$ is a symmetric section of a parabola in ${\cal P}_3$, with endpoints $(1,0,0)$ and $(0,0,1)$ and stationary point $(\frac{1}{4},\frac{1}{2},\frac{1}{4})$.

Because ${\cal S}_{k,n}$ is finite, the uniform distribution $\mathrm{Unif}(k,n)$ on ${\cal S}_{k,n}$ exists and is given by the constant pmf
\begin{equation}\label{unifa}
f_\mathrm{unif}({\bf x})=\ts1/{n+k-1\choose k-1},\quad {\bf x}\in{\cal S}_{k,n}.
\end{equation}
This uniform pmf also can be viewed as a vector in ${\cal P}_{n+k-1\choose k-1}$:
\begin{align}
f_\mathrm{unif}&=\ts{\bf 1}/{n+k-1\choose k-1},\label{boldfunif}\\
{\bf 1}:&=(1,\dots,1)\in\mathbb{R}^{n+k-1\choose k-1}.\nonumber
\end{align}
It is important to note that $f_\mathrm{unif}\notin{\cal M}(k,n)$\footnote{\label{separated}Suppose to the contrary that $f_\mathrm{unif}=f_{\bf p}$ for some ${\bf p}$, that is,
\begin{equation*}
{n\choose {\bf x}}\prod\nolimits_{j=1}^k p_j^{x_j} = \mathrm{constant}\ \forall {\bf x}\in{\cal S}_{k,n}.
\end{equation*}
Take $(x_1,\dots,x_k)=(0,\dots,0,n,0,\dots,0)$ to see that $p_1^n=\cdots=p_k^n$, so $p_1=\cdots=p_k$. Then take $(x_1,\dots,x_k)=(n-1,1,0,\dots,0)$ to see that $p_1^n=np_1^{n-1}p_2$, so $p_1=np_2\ne p_2$, a contradiction.} 
but $f_\mathrm{unif}$ lies in the interior of the convex hull of ${\cal M}(k,n)$. The latter follows from the well-known fact that if $\nu$ is the uniform distribution on ${\cal P}_k$ then
\begin{equation}\label{unifintegral}
{n\choose {\bf x}}\int_{{\cal P}_k}\left(\prod\nolimits_{j=1}^k p_j^{x_j}\right)d\nu({\bf p})=\ts1/{n+k-1\choose k-1}.
\end{equation}


\nid{\bf 3. The first testing problem.}

\nid Let ${\bf X}\equiv(X_1,\dots,X_k)$ be a rvtr with pmf $f$ on the sample space ${\cal S}_{k,n}$. We shall study two problems. The first is that of testing the simple hypothesis
\begin{equation}\label{PST3}
H_{\bf p}:f=f_{\bf p}\quad\mathrm{vs.}\quad H_\mathrm{unif}:f= f_\mathrm{unif}
\end{equation}
for any fixed ${\bf p}\in{\cal P}_k$. The classical LRT for this problem rejects $H_{\bf p}$ for large values of $f_\mathrm{unif}({\bf X})/f_{\bf p}({\bf X})$, equivalently, for small values of $f_{\bf p}({\bf X})$, so is equivalent to the PST for $H_{\bf p}$. If ${\bf X}={\bf x}_0$ is observed, the attained $p$ value is
\begin{equation}\label{PST4}
\pi_{\bf p}({\bf x_0}):=\P_{\bf p}[f_{\bf p}({\bf X})\le f_{\bf p}({\bf x}_0)].
\end{equation}

Determination of $\pi_{\bf p}({\bf x_0})$ is a challenging computational question that is not addressed here. However, Proposition 3.1 and Conjecture 3.3 (if true) indicate that the LRT $\equiv$ PST for \eqref{PST3} when ${\bf p}\ne{\bf p}_\mathrm{ecp}$ is more sensitive, i.e., more powerful (see Remark 3.2) than the LRT $\equiv$ PST when ${\bf p}={\bf p}_\mathrm{ecp}$, the equal-cell-probability (ecp) case,  where
\begin{align}
{\bf p}_\mathrm{ecp}&=\ts(\frac{1}{k},\dots,\frac{1}{k})\in{\cal P}_k,\label{ecp2}\\
f_\mathrm{ecp}({\bf x})&={n\choose{\bf x}}\frac{1}{k^n},\ \ {\bf x}\in{\cal S}_{k,n}.\label{fecp6}
\end{align}
This test rejects $H_\mathrm{ecp}\equiv H_{{\bf p}_\mathrm{ecp}}$ in favor of $H_\mathrm{unif}$ for large values of
\begin{equation}\label{optimaltest}
\frac{f_\mathrm{unif}({\bf X})}{f_\mathrm{ecp}({\bf X})}\propto {n\choose{\bf X}}^{-1}\!\!\!\propto\ts\ \prod_{j=1}^k X_j!\ \ .
\end{equation}
Thus the ecp case should provide a floor for the general case.
\vskip4pt

\nid{\bf Proposition 3.1.} Under the KLD criterion, $f_\mathrm{ecp}$
is the closest pmf in ${\cal M}(k,n)$ to $f_\mathrm{unif}$.
That is,
\begin{equation}\label{KL1}
\E_\mathrm{unif}\left[\log\frac{f_\mathrm{unif}({\bf X})}{f_\mathbf{p}({\bf X})}\right]>\E_\mathrm{unif}\left[\log\frac{f_\mathrm{unif}({\bf X})}{f_\mathrm{ecp}({\bf X})}\right]\ \ \forall {\bf p}\in{\cal P}_k\stms\{{\bf p}_\mathrm{ecp}\}.
\end{equation}
The left side and right side of \eqref{KL1} are the (positive) KLDs from $f_{\bf p}$ to $f_\mathrm{unif}$ and from $f_\mathrm{ecp}$ to $f_\mathrm{unif}$, respectively.
\vskip2pt

\nid{\bf Proof.} Inequality \eqref{KL1} is derived as follows. From \eqref{pmfa} and \eqref{fecp6},
\begin{align*}
&\E_\mathrm{unif}\left[\log \frac{f_\mathrm{unif}({\bf X})}{f_{\bf p}({\bf X})}\right]-\E_\mathrm{unif}\left[\log \frac{f_\mathrm{unif}({\bf X})}{f_\mathrm{ecp}({\bf X})}\right]\\
=\ &\E_\mathrm{unif}\left[\log \frac{f_\mathrm{ecp}({\bf X})}{f_{\bf p}({\bf X})}\right]\\
=\ &-\E_\mathrm{unif}\left[\log k^n\prod_{j=1}^k p_j^{X_j}\right]\\
=\ &-\left[n\log k+\sum_{j=1}^k \E_\mathrm{unif}(X_j)\log p_j\right]\\
=\ &-n\left[\log k+\frac{1}{k}\sum_{j=1}^k \log p_j\right]>0
\end{align*}
by the symmetry of $f_\mathrm{unif}$ and the strict convexity of $\log x$.\hfill$\square$
\vskip6pt

\nid{\bf Remark 3.2.} Chernoff (1956) showed that the logarithm of the Type 2 error probability of the LRT\ $\equiv$ PST is asymptotically proportional to $-$KLD.\hfill$\square$ 
\vskip6pt

Now consider the EPV criterion for the PSTs. For a discussion of the role of the EPV in hypothesis testing, see Sackrowitz and Samuel-Cahn (1999). Abbreviate ${\cal S}_{k,n}$ by ${\cal S}$, let ${\bf Y}\sim f_\mathrm{unif}\ \ \text{on}\ {\cal S}$, and consider the $p$-values $\pi_{\bf p}$ and $\pi_\mathrm{unif}\equiv\pi_{{\bf p}_\mathrm{unif}}$ given by \eqref{PST4}.
\vskip6pt

\nid{\bf Conjecture 3.3:} Under the EPV criterion, $f_\mathrm{ecp}$
is the closest pmf in ${\cal M}(k,n)$ to $f_\mathrm{unif}$. That is, the EPV $\E_\mathrm{unif}[\pi_{\bf p}({\bf Y})]$ is maximized when $\mathbf{p}={\bf p}_\mathrm{ecp}$, i.e., 
 \begin{equation}\label{EpiY}
 \hskip30pt\E_\mathrm{unif}[\pi_\mathrm{ecp}({\bf Y})]>\E_\mathrm{unif}[\pi_{\bf p}({\bf Y})]\ \ \forall {\bf p}\in{\cal P}_k\stms\{{\bf p}_\mathrm{ecp}\}.
 \hskip70pt\square
\end{equation}

From \eqref{PST4}, under the uniform alternative $H_\mathrm{unif}$, the EPV of the LRT $\equiv$ PST for \eqref{PST3} with a general ${\bf p}$  is given by
\begin{align*}
 \E_\mathrm{unif}[\pi_{\bf p}({\bf Y})]
 &=\E_\mathrm{unif}\{\P_{\bf p}\left[ f_{\bf p}({\bf X})\le f_{\bf p}({\bf Y})\,|\,{\bf Y}\right]\}\\
  &=\ts\P_{{\bf p},\mathrm{unif}}\left[ f_{\bf p}({\bf X})\le f_{\bf p}({\bf Y})\right]\\
  &=\ts{n+k-1\choose k-1}^{-1}\sum_{\left\{({\bf x},{\bf y})\in\times{\cal S}\big| f_{\bf p}({\bf x})\le f_{\bf p}({\bf y})  \right \}} f_{\bf p}({\bf x}),
\end{align*}
while under $H_\mathrm{unif}$, the EPV of the LRT $\equiv$ PST for \eqref{PST3} with ${\bf p}={\bf p}_\mathrm{unif}$  is
\begin{align*}
   \E_\mathrm{unif}[\pi_\mathrm{ecp}({\bf Y})]
  &=\ts{n+k-1\choose k-1}^{-1}\sum_{\left\{({\bf x},{\bf y})\in{\cal S}\times{\cal S}\big| f_\mathrm{ecp}({\bf x}) \le f_\mathrm{ecp}({\bf y})\right \}} f_\mathrm{ecp}({\bf x}).
\end{align*}
Thus the conjectured inequality \eqref{EpiY} holds iff
\begin{align}
\sum\limits_{\left\{({\bf x},{\bf y})\in{\cal S}\times{\cal S}\big|{n\choose{\bf x}}\le{n\choose{\bf y}} \right \}} f_\mathrm{ecp}({\bf x})\ \ 
>\sum\limits_{\left\{({\bf x},{\bf y})\in{\cal S}\times{\cal S}\big| f_{\bf p}({\bf x})\le f_{\bf p}({\bf y})  \right \}} f_{\bf p}({\bf x}).\label{ineq5}
\end{align}

By conditioning on ${\bf x}$, inequality \eqref{ineq5} can be rewritten as
\begin{align}
&\sum\limits_{\left\{{\bf x}\in{\cal S}\right \}} \left|\left\{{\bf y}\in{\cal S}\Big|{n\choose{\bf x}}\le{n\choose{\bf y}} \right\}\right|f_\mathrm{ecp}({\bf x})\nonumber\\
>&\sum\limits_{\left\{{\bf x}\in{\cal S}\right \}} \left|\left\{{\bf y}\in{\cal S}\Big|f_{\bf p}({\bf x})\le f_{\bf p}({\bf y})\right\}\right|f_{\bf p}({\bf x}),\label{ineq6}
\end{align}
where $|A|$ denotes the number of elements of $A$. Verification of \eqref{ineq6} has proven to be elusive, although it is straightforward for the extreme case where ${\bf p}$ is any of the $n$ permutations of ${\bf p}_0\equiv(1,0,\dots,0)$:
\vskip6pt

\nid{\bf Proposition 3.4.} Inequality \eqref{ineq6} holds when ${\bf p}$ is any permutation of ${\bf p}_0$.
\vskip2pt

\nid{\bf Proof.} Because ${\cal S}$ is invariant under permutations of ${\bf x}$, the right-hand side of \eqref{ineq6} is invariant under permutations of ${\bf p}$, so it suffices to consider ${\bf p}={\bf p}_0$. Because $f_{{\bf p}_0}({\bf x})=1(0)$ if ${\bf x}=(\ne)\,{\bf x}_0:=(n,0,\dots,0)$, the right-hand side of \eqref{ineq6} $=1$ when ${\bf p}={\bf p}_0$. But the left-hand side $>1$, since 
\begin{align}
\ts\left|\left\{{\bf y}\in{\cal S}\Big|{n\choose{\bf x}}\le{n\choose{\bf y}} \right\}\right|&\ge1\quad\forall\,{\bf x}\in{\cal S},\nonumber
\end{align}
while ${n\choose{\bf x}_0}=1$ so
\begin{align}
\hskip110pt\ts\left|\left\{{\bf y}\in{\cal S}\Big|{n\choose{\bf x}_0}\le{n\choose{\bf y}} \right\}\right|&=|{\cal S}|>1.\nonumber\hskip100pt\square
\end{align}

After Proposition 3.4, the next simplest case of Conjecture 3.3 occurs when ${\bf p}$ has exactly two positive components, which by symmetry can be taken to be $p_1,p_2>0$, with $p_3=\cdots=p_k=0$. This reduces the multinomial family to the binomial family, which will be examined in Sections 5-7, providing further evidence for the validity of Conjecture 3.3 in general.
\vskip6pt
  
\nid{\bf Remark 3.5.} Since we are testing $f= f_{\bf p}\ \mathrm{vs.}\ f= f_\mathrm{unif}$ (cf. \eqref{PST3}), it would also seem to be of interest to compare the expected $p$-values under the null hypotheses $f_\mathrm{ecp}$ and $f_{\bf p}$, that is to compare $ \E_\mathrm{ecp}[\pi_\mathrm{ecp}({\bf Y})]$ and $\E_{\bf p}[\pi_{\bf p}({\bf Y})]$. However, if ${\bf Y}$ had a continuous distribution then $\pi_\mathrm{ecp}({\bf Y})$ would have the Uniform$(0,1)$ distribution under $f_\mathrm{ecp}$, as would $\pi_{\bf p}({\bf Y})$ under $f_{\bf p}$, hence $\E_\mathrm{ecp}[\pi_\mathrm{ecp}({\bf Y})]=\E_{\bf p}[\pi_{\bf p}({\bf Y})]=\t12$. Thus under the actual discrete distribution of ${\bf Y}$, both these expectations are approximately $\t12$, at least for moderate or large $n$, so this comparison would be uninformative.\hfill$\square$
\vskip6pt

\nid{\bf 4. The second testing problem.}

\nid Now consider the problem of testing the composite null hypothesis
\begin{equation}\label{HvHm}
H_\mathrm{mult}: f\in{\cal M}(k,n)\quad\mathrm{vs.}\quad H_\mathrm{unif}: f=f_\mathrm{unif}.
\end{equation}
This can be viewed as testing goodness-of-fit for the multinomial distribution, but with a specified alternative, namely the uniform distribution.
Because 
$f_\mathrm{unif}\notin{\cal M}(k,n)$ 
while $f_\mathrm{ecp}$ is the pmf in ${\cal M}(k,n)$ closest to $f_\mathrm{unif}$ according to KLD (Proposition 3.1) and possibly EPV (Conjecture 3.3), one might ask if $f_\mathrm{ecp}$ is the least favorable distribution for the testing problem \eqref{HvHm}.\footnote{\label{leastfav}Phrased more properly, is the prior distribution on ${\cal M}(k,n)$ that assigns mass 1 to $f_\mathrm{ecp}$ a least favorable prior distribution for \eqref{HvHm} (see [LR], \S3.8)).}
If this were true then the LRT $\equiv$ PST \eqref{optimaltest} would be the most powerful test of its size for \eqref{HvHm}. By Corollary 3.8.1 of [LR] Lehmann and Romano (2005), $f_\mathrm{ecp}$ would be the least favorable distribution if, for all $c>0$,
\begin{equation*}
\P_{\bf p}[f_\mathrm{ecp}({\bf X})\le c]\le\P_\mathrm{ecp}[f_\mathrm{ecp}({\bf X})\le c]\ \ \ \forall {\bf p}\not\in{\cal P}_k.
\end{equation*}
However, exactly the opposite is true:
\vskip4pt

\nid{\bf Proposition 4.1.} $\P_{\bf p}[f_\mathrm{ecp}({\bf X})\le c]$ is a Schur-convex\footnote{\label{majorization}Refer to [MOA] Marshall, Olkin, Arnold (2011) for definitions of Schur-convexity, Schur-concavity, majorization, and $T$-transforms.} function of ${\bf p}$. Thus
\begin{equation}\label{probineq}
\P_{\bf p}[f_\mathrm{ecp}({\bf X})\le c]\ge\P_\mathrm{ecp}[f_\mathrm{ecp}({\bf X})\le c]\ \ \ \forall {\bf p}\notin{\cal P}_k.
\end{equation}
\vskip2pt
\nid{\bf Proof.} From \eqref{fecp6},
\begin{align}
\P_{\bf p}[f_\mathrm{ecp}({\bf X})\le c]&=\ts\P_{\bf p}[{n\choose{\bf X}} \le c']=\P_{\bf p}[\sum_{j=1}^k\log (X_j!)\ge c'']\label{Schur1}
\end{align}
for some constants $c'$ and $c''$. Since $x!=\Gam(x+1)$ and the gamma function is log convex, $\sum_{j=1}^k\log(x_j!)$ is convex and symmetric in ${\bf x}$, so is a Schur-convex function. This implies that the indicator function of the set $\{{\bf x}\,|\,\sum_{j=1}^k\log(x_j!)< c''\}$ is Schur-concave in ${\bf x}$, hence its expected value under $\P_{\bf p}$ is Schur-concave in ${\bf p}$ by the Proposition in Example 2 of [R] Rinott (1973) (see [MOA] Proposition 11.E.11.). By \eqref{Schur1} this expected value is
\begin{align*}
\P_{\bf p}[f_\mathrm{ecp}({\bf X})> c]=1-\P_{\bf p}[f_\mathrm{ecp}({\bf X})\le c],
\end{align*}
so the first assertion holds. Because every ${\bf p}$ majorizes ${\bf p}_\mathrm{ecp}$, \eqref{probineq} follows.\hfill$\square$
\vskip6pt

In view of Proposition 4.1, the LRT $\equiv$ PST \eqref{optimaltest} is inappropriate for \eqref{HvHm}. Instead we propose the actual LRT for \eqref{HvHm}, which rejects $H_\mathrm{mult}$ in favor of $H_\mathrm{unif}$ for {\it large} values of
\begin{equation}\label{LRT}
\frac{f_\mathrm{unif}({\bf X})}{\sup\limits_{\mathbf{p}\in{\cal P}_k}f_\mathbf{p}({\bf X})}\propto\prod_{j=1}^k \frac{X_j!}{X_j^{X_j}}=:L({\bf X}).
\end{equation} 
By Proposition 4.2, this test has the desirable property that its power function is a Schur-concave function of ${\bf p}$, hence attains its maximum over the null hypothesis $H_\mathrm{mult}$ at ${\bf p}={\bf p}_\mathrm{ecp}$, in conformity with the minimum property of the KLD in \eqref{KL1}. This implies that the Type 2 error probability attains its maximum over $H_\mathrm{mult}$ at ${\bf p}$, hence   this LRT is unbiased for \eqref{HvHm}.
\vskip4pt

\nid{\bf Proposition 4.2.} $\P_{\bf p}[L({\bf X})\ge c]$ is a Schur-concave function of ${\bf p}\in{\cal P}_k$.
\vskip2pt
\nid{\bf Proof.} We will show that $L({\bf x})$ is Schur-concave, hence the indicator function $I_c({\bf x})$ of the set $\{{\bf x}\,|\,L({\bf x})\ge c\}$ is also Schur-concave. Then by the Proposition in Example 2 of [R], $\E_{\bf p}[I_c({\bf x})]\equiv\P_{\bf p}[L({\bf X})\ge c]$ is Schur-concave in ${\bf p}$.

To show that $L({\bf x})$ is Schur-concave, we must show that $L({\bf y})\le L({\bf x})$ whenever ${\bf y}$ majorizes ${\bf x}$. By a classical result of Muirhead (cf. [MOA] Lemma 3.1), ${\bf x}$ can be obtained from ${\bf y}$ by a finite number of linear $T${\it -transforms}, which act on pairs of the coordinates of ${\bf y}$ by moving both members of the pair toward their average. Thus by the symmetry of $L$, it suffices to show
\begin{equation}\label{xy}
 L(x_1-1,x_2+1,x_3,\dots,x_k)\le L(x_1,x_2,x_3,\dots,x_k)
\end{equation}
when $x_1\le x_2$. But this is equivalent to each of the inequalities
\begin{align}
\frac{(x_1-1)!(x_2+1)!}{(x_1-1)^{x_1-1}(x_2+1)^{x_2+1}}&\le\frac{x_1!x_2!}{x_1^{x_1}x_2^{x_2}},\nonumber\\
\ts\left(1+\frac{1}{x_1-1}\right)^{x_1-1}&\ts\le\left(1+\frac{1}{x_2}\right)^{x_2},\label{ineq88}
\end{align}
and \eqref{ineq88} holds since $\left(1+\frac{1}{t}\right)^{t}$ increases for $t>0$ and $x_1-1<x_1\le x_2$.\hfill$\square$
\vskip4pt

\nid{\bf Remark 4.3.} If ${\bf X}={\bf x}_0$ is observed, the attained $p$-value of the LRT \eqref{LRT} can be expressed in a tractable form:
\begin{align}
\pi_L({\bf x}_0):&=\sup_{{\bf p}\in{\cal P}_k}\P_{\bf p}[L({\bf X})\ge L({\bf x}_0)]\nonumber\\
 &=\P_\text{ecp}[L({\bf X})\ge L({\bf x}_0)].\label{pLRT}
\end{align}
Because ${\bf p}$ majorizes ${\bf p}_\text{ecp}$\ $\forall{\bf p}\in{\cal P}_k$, \eqref{pLRT} follows from Proposition 4.2. Evaluation of  \eqref{pLRT}, based on $f_\mathrm{ecp}({\bf x})$,  is not addressed here.\hfill$\square$ 
\vskip4pt

\nid{\bf Remark 4.4.} Evaluation of \eqref{pLRT} is easy in the binomial case $k=2$, where ${\bf x}=(x_1,n-x_1)$, ${\bf x}_0=(x_{10},n-x_{10})$, and 
\begin{equation}\label{LBin}
L(x_1,n-x_1)=\frac{x_1!(n-x_1)!}{x_1^{x_1}(n-x_1)^{n-x_1}},\quad x_1=0,\dots,n.
\end{equation}
Here the Schur-concavity of $L({\bf x})$, shown in the proof of Proposition 4.2, implies that  $L(x_1,n-x_1)$ is unimodal and symmetric about $x_1=n/2$. Therefore 
\begin{equation}\label{ineq86}
\P_\text{ecp}[L({\bf X})\ge L({\bf x}_0)]=\ts\P_\text{ecp}[\,|X_1-\frac{n}{2}|\le|x_{10}-\frac{n}{2}|\,].
\end{equation}
Because $X_1\sim\text{Binomial}(n,\t12)$ in the ecp case, this is readily evaluated.\hfill$\square$
\vskip6pt

\nid{\bf 5. The  binomial case ($k=2$).} 

\nid In Sections 5-7 we examine the validity of Conjecture 3.3 for the first testing problem \eqref{PST3} in the apparently simple but still challenging binomial case, which hopefully might suggest an approach to the general multinomial case.

Denote the Binomial$(n,p)$ pmf by
\begin{equation}\label{pmfB}
\ts f_p(x)={n\choose x}p^x(1-p)^{n-x},\quad x=0,1\dots,n,\quad0\le p\le1.
\end{equation}
 When $k=2$, the inequality \eqref{ineq6} in Conjecture 3.3 reduces to
  \begin{align}
&\sum_{x=0}^n q_{\frac{1}{2}}(x)f_{\frac{1}{2}}(x)\ge\sum_{x=0}^n  q_p(x)f_p(x),\label{ineq7x}
\end{align}
where 
\begin{equation}\label{rpx}
q_p(x)=\left|\left\{y\,\Big|\,0\le y\le n,\,f_p(x)\le f_p(y)\right\}\right|.
\end{equation}
Note that $1\le q_p(x)\le n+1$.  I believe that \eqref{ineq7x} can be sharpened as follows:
\vskip4pt


 \nid{\bf Conjecture 5.1.} If $\t12<p\le1$ then
 \begin{align}
\hskip80pt\sum_{x=0}^n q_{\frac{1}{2}}(x)f_{\frac{1}{2}}(x)>\sum_{x=0}^n  q_p(x)f_p(x). \hskip100pt\square\label{ineq7}
\end{align}
\vskip4pt

 It is easy to see that
\begin{equation}\label{qabs}
q_{\frac{1}{2}}(x)=1+|n-2x|,
\end{equation}
so \eqref{ineq7} can be written as
 \begin{align}
&1+\sum_{x=0}^n |n-2x|f_{\frac{1}{2}}(x)>\sum_{x=0}^n  q_p(x)f_p(x).\label{ineq6.1}
\end{align}
Because $f_p(x)=f_{1-p}(n-x)$ and $q_p(x)=q_{1-p}(n-x)$, only the case $\t12\le p\le1$ need be considered. The difficulty in verifying \eqref{ineq6.1} stems mainly from the difficulty in determining $q_p(x)$ as $p$ increases from $\t12$ to 1. Two special cases are somewhat amenable (proved in Appendix 1):
\vskip8pt


\nid{\bf Proposition 5.2.} For the binomial case, Conjecture 5.1 holds when
\vskip2pt
\nid (i) $\frac{n}{n+1}\le p\le1$ (i.e., $p$ near the extreme case $p=1$);
\vskip2pt
\nid (ii) $\t12< p<\t12+\del_n$ for some sufficiently small $\del_n\le\eps_n$ (i.e., $p$ near the ecp case $p=\frac{1}{2}$), where
\begin{align}
\hskip50pt \eps_n&=\frac{\c a_n-1}{2[\c a_n+1]};\label{epsn}\\
\c a_n&=\begin{cases}\min\left\{a_n(x)\,\Big|\,x={\ts\frac{n+3}{2}},\dots,n\right\},&\quad n\ \text{odd},\\\min\left\{a_n(x)\,\Big|\,x={\ts\frac{n}{2}+1},\dots,n\right\},&\quad n\ \text{even};\end{cases}\label{epsn1}\\
a_n(x)&=\left(\frac{x}{n-x+1}\right)^{\frac{1}{2x-n-1}}>1.\label{an}\hskip147pt\square
\end{align}

\vskip4pt

\vskip8pt

\nid{\bf 6. Three binomial examples.} 

\nid In the binomial case the conjectured inequality \eqref{ineq7} holds for the first nontrivial cases $n=3,4,5$. To verify this, first determine $q_p(x)$ as in Tables 1, 2, 3 below, then determine the $n$-th degree polynomial $\sum_{x=0}^n  q_p(x)f_p(x)$. For $n=3$ there are six such polynomials, each of degree 3 or less, one for each of the six sub-ranges of $p$ in Table 1. For $n=4$ ($n=5$) there are ten (fourteen) sub-ranges in Table 2 (Table 3) hence ten (fourteen) such polynomials, each of degree 4 (5) or less but not shown. It is straightforward to verify numerically that \eqref{ineq7} holds for each of these polynomials over their respective sub-ranges of $p$ in Tables 1, 2, 3. (However, \eqref{ineq7} usually does not extend from the sub-range to the entire range $[\t12,1]$). 

These three examples begin to illustrate   the complexity of determining the values of $q_p(x)$, but perhaps suggest a pattern which might be extended for general $n$. However, caution must be exercised, as demonstrated now. 

In Tables 1,2,3, as $p$ begins to increase above $\t12$ the values of $q_p(x)$ remain unchanged for $x\le n/2$ and decrease by 1 for $x>n/2$. As $p$ continues to increase, the values of $q_p(x)$ change in pairs: in Table 1, $n=3$ is odd and the first pair to change is $(q_p(1),q_p(3))$, from (2,3) to (3,3), to (3,2); in Table 3, $n=5$ is odd and the first pair to change is $(q_p(2),q_p(4))$, again from (2,3) to (3,3), to (3,2). This suggests that for odd $n$, the first pair to change is always $(q_p(\frac{n-1}{2}),q_p(\frac{n+3}{2}))$. However, this is not the case:
\begin{align*}
q_p(\frac{n-1}{2})=q_p(\frac{n+3}{2})&\iff f_p(\frac{n-1}{2})=f_p(\frac{n+3}{2})\\
&\iff t=\left(\frac{n+3}{n-1}\right)^{\frac{1}{2}},
\end{align*}
where $t=\frac{p}{1-p}$, while
\begin{align*}
q_p(1)=q_p(n)&\iff f_p(1)=f_p(n)\\
&\iff t=n^{\frac{1}{n-1}},
\end{align*}
hence the pair $(q_p(\frac{n-1}{2}),q_p(\frac{n+3}{2}))$ changes before the pair $(q_p(1),q_p(n))$ iff 
\begin{equation*}
\left(\frac{n+3}{n-1}\right)^{\frac{1}{2}}<n^{\frac{1}{n-1}}.
\end{equation*}
This holds for $n\le45$ but fails for $n\ge47$.

\begin{table}[h!]
\centering
\begin{center}
\begin{tabular}{ |c|c|c|c|c|c|c|c|c| }
 \hline
 $t\equiv\frac{p}{1-p}$&$p$&$q_p(0)$&$q_p(1)$&$q_p(2)$&$q_p(3)$&$\sum\limits_{x=0}^3  q_p(x)f_p(x)$\\
 \hline
 \hline
  $1$&$\t12\equiv\frac{1}{1+1}$&4&2&2&4&$4-6p+6p^2$\\
    \hline
 $(1,3^{\frac{1}{2}})$&$(\frac{1}{1+1},\frac{3^{\frac{1}{2}}}{3^{\frac{1}{2}}+1})$&4&2&1&3&$4-6p+3p^2+2xp^3$\\
 \hline
 $3^{\frac{1}{2}}$&$\frac{3^{\frac{1}{2}}}{3^{\frac{1}{2}}+1}$&4&3&1&3&$4-3p-3p^2+5p^3$\\
 \hline
$(3^{\frac{1}{2}},3)$&$(\frac{3^{\frac{1}{2}}}{3^{\frac{1}{2}}+1},\frac{3}{3+1})$&4&3&1&2&$4-3p-3p^2+4p^3$\\
 \hline
 $3$&$\frac{3}{3+1}$&4&3&2&2&$4-3p+p^3$\\
 \hline
 $(3,\infty]$&$(\frac{3}{3+1},1]\ \ $&4&3&2&1&$4-3p$\\
 \hline
\end{tabular}
\end{center}
\caption{Determination of $q_p(x)$ for $n=3$.}
\vskip24pt
\label{table:1}
\end{table}

By contrast, when $n$ is even, the unimodality of $f_p(0),\dots,f_p(n)$ and the symmetry of $f_{\frac{1}{2}}(x)$ about $x=n/2$ implies that the first pair to change must be one of the pairs $(q_p(n-x+1),q_p(x)$, $x=\frac{n}{2}+1,\dots,n$. For such pairs, the change occurs when
\begin{align*}
q_p(n-x+1)=q_p(x)&\iff f_p(n-x+1)=f_p(x)\\
&\iff t=\left(\frac{x}{n-x+1}\right)^{\frac{1}{2x-n-1}}=:Q(x).
\end{align*}
However,
\begin{align*}
\log Q(x)&=\frac{\log x-\log(n-x+1)}{x-(n-x+1)}\\
 &=\frac{1}{2x-n-1}\sum_{m=n-x+2}^x d_m,\\
 d_m:&=\log(m)-\log(m-1)
\end{align*}
and $d_m$ is strictly convex in $m$, hence $\log Q(x)$ is strictly increasing in $x$ for $x=\frac{n}{2}+1,\dots,n$ by Lemma 6.1 below (proved in Appendix 1). This implies that for even $n$, $(q_p(\frac{n}{2}),q_p(\frac{n}{2}+1))$ is always the first pair to change.
\vskip6pt

\begin{table}[h!]
\centering
\begin{center}
\begin{tabular}{ |c|c|c|c|c|c|c|c|c|c| }
 \hline
 $t\equiv\frac{p}{1-p}$&$p$&$q_p(0)$&$q_p(1)$&$q_p(2)$&$q_p(3)$&$q_p(4)$\\
 \hline
 \hline
  $1$&$\t12\equiv\frac{1}{1+1}$&5&3&1&3&5\\
    \hline
 $(1,\frac{3}{2})$&$(\frac{1}{1+1},\frac{\frac{3}{2}}{\frac{3}{2}+1})$&5&3&1&2&4\\
 \hline
 $\frac{3}{2}$&$\frac{\frac{3}{2}}{\frac{3}{2}+1}$&5&3&2&2&4\\
 \hline
$(\frac{3}{2},4^{\frac{1}{3}})$&$(\frac{\frac{3}{2}}{\frac{3}{2}+1},\frac{4^{\frac{1}{3}}}{4^{\frac{1}{3}}+1})$&5&3&2&1&4\\
 \hline
 $4^{\frac{1}{3}}$&$\frac{4^{\frac{1}{3}}}{4^{\frac{1}{3}}+1}$&5&4&2&1&4\\
 \hline
 $(4^{\frac{1}{3}},6^{\frac{1}{2}})$&$(\frac{4^{\frac{1}{3}}}{4^{\frac{1}{3}}+1},\frac{6^{\frac{1}{2}}}{6^{\frac{1}{2}}+1})$&5&4&2&1&3\\
 \hline
 $6^{\frac{1}{2}}$& $\frac{6^{\frac{1}{2}}}{6^{\frac{1}{2}}+1}$&5&4&3&1&3\\
 \hline
 $(6^{\frac{1}{2}},4)$&$(\frac{6^{\frac{1}{2}}}{6^{\frac{1}{2}}+1},\frac{4}{4+1})$&5&4&3&1&2\\
 \hline
 $4$&$\frac{4}{4+1}$&5&4&3&2&2\\
 \hline
$(4,\infty]$&$(\frac{4}{4+1},1]\ \ $&5&4&3&2&1\\
 \hline
\end{tabular}
\end{center}
\caption{Determination of $q_p(x)$ for $n=4$.}
\label{table:2}
\end{table}

\newpage

\nid{\bf Lemma 6.1.} For $n\ge4$ let $d_1,\dots,d_n$
be a strictly convex sequence. For $\frac{n}{2}+1\le x\le n$, define 
\begin{equation*}
D_x=\frac{1}{2x-n-1}\sum_{m=n-x+2}^x d_m.
\end{equation*}
Then $D_x$ is strictly increasing in $x$.\hfill$\square$
\vskip6pt

\begin{table}[h!]
\centering
\begin{center}
\begin{tabular}{ |c|c|c|c|c|c|c|c|c|c|c| }
 \hline
 $t\equiv\frac{p}{1-p}$&$p$&$q_p(0)$&$q_p(1)$&$q_p(2)$&$q_p(3)$&$q_p(4)$&$q_p(5)$\\
 \hline
 \hline
  $1$&$\t12\equiv\frac{1}{1+1}$&  6&4&2&2&4&6\\
    \hline
 $(1,2^{\frac{1}{2}})$&$(\frac{1}{1+1},\frac{2^{\frac{1}{2}}}{2^{\frac{1}{2}}+1})$&6&4&2&1&3&5\\
 \hline
 $2^{\frac{1}{2}}$&$\frac{2^{\frac{1}{2}}}{2^{\frac{1}{2}}+1}$&6&4&3&1&3&5\\
 \hline
$(2^{\frac{1}{2}},5^{\frac{1}{4}})$&$(\frac{2^{\frac{1}{2}}}{2^{\frac{1}{2}}+1},\frac{5^{\frac{1}{4}}}{5^{\frac{1}{4}}+1})$&6&4&3&1&2&5\\
 \hline
 $5^{\frac{1}{4}}$&$\frac{5^{\frac{1}{4}}}{5^{\frac{1}{4}}+1}$&6&5&3&1&2&5\\
 \hline
 $(5^{\frac{1}{4}},2)$&$(\frac{5^{\frac{1}{4}}}{5^{\frac{1}{4}}+1},\frac{2}{2+1})$&6&5&3&1&2&4\\
 \hline
 $2$& $\frac{2}{2+1}$&6&5&3&2&2&4\\
 \hline
 $(2,10^{\frac{1}{3}})$&$(\frac{2}{2+1},\frac{10^{\frac{1}{3}}}{10^{\frac{1}{3}}+1})$&6&5&3&2&1&4\\
 \hline
 $10^{\frac{1}{3}}$&$\frac{10^{\frac{1}{3}}}{10^{\frac{1}{3}}+1}$&6&5&4&2&1&4\\
 \hline
$(10^{\frac{1}{3}},10^{\frac{1}{2}})$&$(\frac{10^{\frac{1}{3}}}{10^{\frac{1}{3}}+1},\frac{10^{\frac{1}{2}}}{10^{\frac{1}{2}}+1})$&6&5&4&2&1&3\\
 \hline
 $10^{\frac{1}{2}}$&$\frac{10^{\frac{1}{2}}}{10^{\frac{1}{2}}+1}$&6&5&4&3&1&3\\
 \hline
$(10^{\frac{1}{2}},5)$& $(\frac{10^{\frac{1}{2}}}{10^{\frac{1}{2}}+1},\frac{5}{5+1})$&6&5&4&3&1&2\\
 \hline
 $5$& $\frac{5}{5+1}$&6&5&4&3&2&2\\
 \hline
 $(5,\infty]$&  $(\frac{5}{5+1},1]\ \ \ $&6&5&4&3&2&1\\
 \hline
\end{tabular}
\end{center}
\caption{Determination of $q_p(x)$ for $n=5$.}
\label{table:3}
\end{table}

\nid{\bf 7. The ordered binomial distribution.} 

\nid There is an interesting relation between the binomial Conjecture 5.1 and what I shall call the {\it ordered binomial distribution (OBD)}, which appears to be new and of interest in its own right. 

For $n=1,2,\dots$ and $0\le p\le1$, let $X\equiv X_{n,p}$ be a random variable having the Binomial$(n,p)$ distribution, with pmf $f_p(x)$ given in \eqref{pmfB}. Rearrange the $n+1$ probabilities $f_p(x)$ in ascending order to obtain
\begin{equation}\label{pmfBO}
\~f_p(0)\le\~f_p(1)\le\cdots\le\~f_p(n-1)\le\~f_p(n).
\end{equation}
The relation between $\~f_p(\cdot)$ and $f_p(\cdot)$ is given by
\begin{equation}\label{fprelation}
\~f_p(r_p(x))=f_p(x),\quad x=0,\dots,n,
\end{equation}
where
\begin{align}
r_p(x)&=\left|\left\{y\,\Big|\,0\le y\le n,\,f_p(y)< f_p(x)\right\}\right|\nonumber\\
&=n+1-q_p(x)\label{rqrel}
\end{align}
and $|A|$ denotes the number of elements of $A$. Here $r_p(x)+1$ is the rank of $f_p(x)$ among $f_p(0),\dots,f_p(n)$, where, in the case of a tie $f_p(x_1)=f_p(x_2)$, the lower rank is assigned to both. Clearly $0\le r_p(x)\le n$.
\vskip6pt

\nid{\bf Definition 7.1.} The {\it ordered binomial distribution (OBD)} is the distribution of the random variable (rv) $\~X\equiv\~X_{n,p}$, defined by 
\begin{equation}\label{pmfBOD}
\ts \P[\~X=x]=\~f_p(x),\qquad x=0,1\dots,n.
\end{equation}
Clearly $\~X_{n,1-p}\eqd\~X_{n,p}$ because $f_p(x)=f_{1-p}(n-x)$, so it suffices to study the OBD only for $\t12\le p\le1$. 
By \eqref{rqrel}, Conjecture 5.1 can be restated equivalently as follows:
\vskip4pt

 \nid{\bf Conjecture 7.2.} If $\t12<p\le1$ then
 \begin{align}
\hskip81pt\sum_{x=0}^n  r_p(x)f_p(x)>\sum_{x=0}^n r_{\frac{1}{2}}(x)f_{\frac{1}{2}}(x).\hskip100pt\square\label{ineq71}
\end{align}

 By \eqref{fprelation}, 
 \begin{align}
\sum_{x=0}^n  r_p(x)f_p(x)
=\sum_{x=0}^n  r_p(x)\~f_p(r_p(x)).\label{ftilderel}
\end{align}
Lemma 7.3 (proved in Appendix 1) shows that the sums in \eqref{ftilderel} are either equal to, or very close to, $\E(\~X_p)$.
\vskip4pt

\nid{\bf Lemma 7.3.} (i) $\sum\limits_{x=0}^n  r_p(x)f_p(x)=\E(\~X_p)$ if no ties occur\footnote{Because $f_p(x)$ is a polynomial in $p$, no ties can occur among $f_p(0),\dots,f_p(n)$ if $p$ is non-algebraic, i.e., transcendental, and almost all real numbers are transcendental.} among $f_p(0),\dots,f_p(n)$. This includes the case $\frac{n}{n+1}<p<1$.

\nid(ii) $\sum\limits_{x=0}^n r_{\frac{1}{2}}(x)f_{\frac{1}{2}}(x)=\E(\~X_{\frac{1}{2}})-\t12$.
\vskip2pt

\nid(iii) $\sum\limits_{x=0}^n  r_p(x)f_p(x)=\E(\~X_p)-\Del_{n,p}$ if $\t12<p\le\frac{n}{n+1}$ and ties occur among among $f_p(0),\dots,f_p(n)$; here $\Del_{n,p}$ is some number such that $0<\Del_{n,p}<\t12$.
\vskip2pt

\nid(iv) $\sum\limits_{x=0}^n  r_1(x)f_1(x)=\E(\~X_1)$.$\hfill\square$
\vskip8pt

By Lemma 7.3, the inequality \eqref{ineq71} in Conjecture 7.2 can be restated equivalently in terms of the OBD rv $\~X_p$: 
\begin{equation}\label{ineq72}
\E(\~X_p)>\begin{cases}\E(\~X_{\frac{1}{2}})-\t12&\text{if}\ \t12<p<1\ \text{and\ no\ ties\ occur},\\
        \E(\~X_{\frac{1}{2}})-\Lam_{n,p}&\text{if}\ \t12<p<1\ \text{and\ ties\ occur},\\
        \E(\~X_{\frac{1}{2}})-\t12&\text{if}\ p=1,\end{cases}
\end{equation}
where  $\Lam_{n,p}=\t12-\Del_{n,p}$, so $0<\Lam_{n,p}<\t12$. By Proposition 5.2, \eqref{ineq72} holds for $p$ near $\t12$ and near 1. The three examples show that \eqref{ineq72} holds for $n=3,4,5$.
\vskip2pt

Lastly, the conjectured inequality \eqref{ineq72} can be strengthened successively:
\vskip4pt

\nid{\bf Conjecture 7.4.} $\E(\~X_p)>\E(\~X_{\frac{1}{2}})$ for $\t12<p\le1$.\hfill$\square$
\vskip4pt

\nid{\bf Conjecture 7.5.} $\E(\~X_p)$ is strictly increasing for $\t12\le p\le1$.\hfill$\square$
\vskip4pt

\nid{\bf Conjecture 7.6:} The OBD rv $\~X_p$ is strictly stochastically increasing in $p$. Equivalently, the OBD probability vector $(\~f_p(0),\dots,\~f_p(n))$ is strictly increasing in the majorization ordering for $\t12\le p\le1$.\hfill$\square$
\vskip8pt

I will continue to investigate these conjectures for the binomial case, hopefully to gain some insight for the general multinomial case. For now I conclude with two miscellaneous remarks about the OBD.
\vskip6pt

\newpage

\nid{\bf Remark 7.7.} Lemma 7.3(ii) provides an explicit expression for $\E(\~X_{\frac{1}{2}})$:
Because $r_{\frac{1}{2}}(x)=n-|n-2x|$ by \eqref{qabs} and \eqref{rqrel},  \eqref{lem1} and Lemma 7.3(ii) 
yield 
\begin{align*}
\hskip40pt\sum_{x=0}^n r_{\frac{1}{2}}(x)f_{\frac{1}{2}}(x)&=n-\frac{1}{2^n}\sum_{x=0}^n|n-2x|{n\choose x}\\
&=\begin{cases}n\left[1-\frac{1}{2^n}{n\choose\frac{n}{2}}\right]&\text{if}\ n\ \text{is\ even},\\n\left[1-\frac{1}{2^{n-1}}{n-1\choose\frac{n-1}{2}}\right]&\text{if}\ n\ \text{is\ odd};
\end{cases}\\
\E(\~X_{\frac{1}{2}})&=\begin{cases}n\left[1-\frac{1}{2^n}{n\choose\frac{n}{2}}\right]+\t12&\text{if}\ n\ \text{is\ even},\\n\left[1-\frac{1}{2^{n-1}}{n-1\choose\frac{n-1}{2}}\right]+\t12&\text{if}\ n\ \text{is\ odd}.\hskip56pt\square
\end{cases}
\end{align*}

\nid{\bf Remark 7.8.} For $\t12\le p\le1$, the minimum OBD probability $\~f_p(0)=(1-p)^n$, the minimum binomial probability. The maximum OBD probability $\~f_p(n)=f_p(\^x_{n,p})$, the maximal $\equiv$ modal binomial probability, where the binomial mode $\^x_{n,p}$ satisfies
\begin{align}
f_p(\^x_{n,p})&\ge\max(f_p(\^x_{n,p}-1),f_p(\^x_{n,p}+1)),\nonumber\\
\text{equivalently},\hskip10pt&\ \ \frac{\^x_{n,p}}{n+1}\le p\le\frac{\^x_{n,p}+1}{n+1}.\label{modeineqF}
\end{align}
Thus $\^x_{n,p}$ is the unique integer in the interval $((n+1)p-1,(n+1)p)$ if $(n+1)p$ is not an integer, while $\^x_{n,p}$ occurs at both $(n+1)p-1$ and $(n+1)p$ if $(n+1)p$ is  an integer. In particular,
\begin{align*}
\hskip12pt\~f_{\t12}(n)&=\begin{cases}{n\choose n/2}/2^n&\hskip23pt\text{if}\ n\ \text{is\ even},\\
    {n\choose (n-1)/2}/2^n={n\choose (n+1)/2}/2^n&\hskip23pt\text{if}\ n\ \text{is\ odd};\end{cases}\\
    \~f_p(n)&={n\choose (n+1)p-1}p^{(n+1)p-1}(1-p)^{n-(n+1)p+1}\\
    &={n\choose (n+1)p}p^{(n+1)p}(1-p)^{n-(n+1)p}\quad\text{if}\ (n+1)p\ \text{is\ an\ integer}.\hskip18pt  \square
\end{align*}
\vskip6pt
\newpage
\nid{\bf 8. Repeated observations.} Now suppose  we observe the random matrix
\begin{align}
{\bf X}=\begin{pmatrix}{\bf X}_1\\ \vdots\\ {\bf X}_r\end{pmatrix}
=\begin{pmatrix}X_{11}&\cdots &X_{1k}\\ \vdots & & \vdots\\ X_{r1}&\cdots& X_{rk}\end{pmatrix},\label{repeated}
\end{align}
where ${\bf X}_1,\dots,{\bf X}_r$ are independent, identically distributed rvtrs with common pmf $f$ on ${\cal S}_{k,n}$. The range of ${\bf X}$ is
\begin{equation*}
{\cal S}_{k,n}^r:={\cal S}_{k,n}\times\cdots\times{\cal S}_{k,n}\quad (r\ \mathrm{times}).
\end{equation*}
Reconsider the testing problems \eqref{PST3} in Section 3 and \eqref{HvHm} in Section 4, repeated here for convenience: test
\begin{align}
H_{\bf p}:&\ f =f_{\bf p}&\mathrm{vs.} &&H_\mathrm{unif}:f=f_\mathrm{unif},&&\label{PST38}\\
H_\mathrm{mult}:&\ f\in{\cal M}(k,n)&\mathrm{vs.}&& H_\mathrm{unif}:f=f_\mathrm{unif},&&\label{HvHm8}
\end{align}
but now based on the repeated observations represented by ${\bf X}$. 

If we define
\begin{equation}\label{Xplus}
{\bf X}_+:=\sum_{i=1}^r{\bf X}_i\equiv (X_{+1},\dots,X_{+k}),
\end{equation} 
then ${\bf X}_+$ is a sufficient statistic for ${\bf X}$ under the multinomial model ${\cal M}(k,n)$ and ${\bf X}_+\sim\text{Mult}_{\bf p}(k,r n)$ when $f =f_{\bf p}$.  Therefore it may seem that the results of Sections 3 and 4 apply directly: just replace $n$ and ${\bf X}$ therein by $rn$  and ${\bf X}_+$ throughout. For example, for problem \eqref{HvHm8} the LRT statistic $L({\bf X})$ in \eqref{LRT} apparently would be replaced by
\begin{equation}\label{LRT1}
\frac{f_\mathrm{unif}({\bf X}_+)}{\sup\limits_{\mathbf{p}\in{\cal P}_k}f_\mathbf{p}({\bf X}_+)}\propto\prod_{j=1}^k \frac{X_{+j}!}{X_{+j}^{X_{+j}}}=:\~L({\bf X}_+).
\end{equation} 
However, this approach is incorrect, both because ${\bf X}_+$ is not uniformly distributed on ${\cal S}_{k,rn}$ under $H_\mathrm{unif}$ and because ${\bf X}_+$ is not a sufficient statistic under the combined model ${\cal M}(k,n)\cup\{f_\mathrm{unif}\}$; information would be lost by considering ${\bf X}+$ alone.
Therefore, tests for \eqref{PST38} and \eqref{HvHm8} must be based on the pmf of ${\bf X}$ itself, whose range is ${\cal S}_{k,n}^r$ not ${\cal S}_{k,rn}$. 

For ${\bf x}\in{\cal S}_{k,n}^r$, the pmfs of ${\bf X}$ under $H_{\bf p}$, $H_\mathrm{ecp}$, and $H_\mathrm{unif}$ respectively, are
\begin{align}
f_{\bf p}({\bf x})&=\prod_{i=1}^r\left[{n\choose {\bf x}_i}\prod_{j=1}^kp_j^{x_{ij}}\right]\nonumber\\
 &=\left[\prod_{i=1}^r{n\choose {\bf x}_i}\right]\left[\prod_{j=1}^kp_j^{x_{+j}}\right],\label{repeated5}\\
f_\mathrm{ecp}({\bf x})&=\frac{1}{k^{rn}}\prod_{i=1}^r{n\choose {\bf x}_i},\label{repeated2}\\
f_\mathrm{unif}({\bf x})&={n+k-1\choose k-1}^{-r}.\label{repeated3}
\end{align}
Because $f_\mathrm{unif}({\bf x})$ is uniform on ${\cal S}_{k,n}^r$, the LRT $\equiv$ PST for \eqref{PST38} based on ${\bf X}$ rejects $H_{\bf p}$ for small values of $f_{\bf p}({\bf X})$. Because the KLD based on ${\bf X}$ is $r\,\times$ the KLD based on a single observation ${\bf X}_i$, Proposition 3.1 remains valid here, again indicating that this test is least sensitive when ${\bf p}={\bf p}_\mathrm{ecp}$. For this case the LRT $\equiv$ PST rejects $H_\mathrm{ecp}$ for large values of
\begin{equation}\label{repeated4}
\prod_{i=1}^r{n\choose {\bf X}_i}^{-1}\propto\ \prod_{i=1}^r\prod_{j=1}^k X_{ij}! \ ,
\end{equation}
which should be compared to \eqref{optimaltest}. We believe that Conjecture 3.3 regarding the EPV criterion is also valid for repeated observations.

Next, the LRT for \eqref{HvHm8} based on ${\bf X}$ rejects $H_\mathrm{mult}$ for large values of 
\begin{equation}\label{LRT9}
\frac{f_\text{unif}({\bf X})}{\sup\limits_{{\bf p}\in{\cal P}_k} f_{\bf p}({\bf X})}\propto\frac{\prod_{i=1}^r\prod_{j=1}^kX_{ij}!}{\prod_{j=1}^kX_{+j}^{X_{+j}}}=:L^*({\bf X}),
\end{equation} 
which should be compared to \eqref{LRT} and \eqref{LRT1}. We conjecture that Proposition 4.2 extends to the repeated observations case as follows:
\vskip4pt

\nid{\bf Conjecture 8.1.} $\P_{\bf p}[L^*({\bf X})\ge c]$ is a Schur-concave function of ${\bf p}\in{\cal P}_k$.\hfill$\square$
\vskip4pt

 If Conjecture 8.1 is true then as in Remark 4.3, if ${\bf X}={\bf x}_0$ is observed, the attained $p$-value of the LRT \eqref{LRT9} can be expressed in a tractable form:
\begin{align}
\pi_{L^*}({\bf x}_0)&=\sup_{{\bf p}\in{\cal P}_k}\P_{\bf p}[L^*({\bf X})\ge L^*({\bf x}_0)]\nonumber\\
 &=\P_\text{ecp}[L^*({\bf X})\ge L^*({\bf x}_0)].\label{pLRT9}
\end{align}
Here \eqref{pLRT9} would follow from Conjecture 8.1 since ${\bf p}$ majorizes ${\bf p}_\text{ecp}$\ $\forall{\bf p}\in{\cal P}_k$. As with \eqref{pLRT}, evaluation of  \eqref{pLRT9} is not addressed here.

Suppose, however, that we attempt to apply the method used for the proof of Proposition 4.2 in order to verify Conjecture 8.1. Write
\begin{align}
\P_{\bf p}[L^*({\bf X})\ge c]&=\E_{\bf p}\{\P[L^*({\bf X})\ge c\,|\,{\bf X}_+]\},\label{PEP}
\end{align}
where the conditional pmf $f_\text{mult}({\bf x}\,|\,{\bf x}_+)$ under the multinomial model ${\cal M}(k,n)$ does not depend on ${\bf p}$:
\begin{align}
f_\text{mult}({\bf x}\,|\,{\bf x}_+)&=\frac{\prod_{i=1}^r{n\choose{\bf x}_i}p_j^{x_{ij}}}{{rn\choose {\bf x}_+}\prod_{j=1}^kp_j^{x_{+j}}}
 =\frac{\prod_{i=1}^r{n\choose{\bf x}_i}}{{rn\choose {\bf x}_+}},\qquad{\bf x}\in{\cal R}({\bf x}_+),\label{condpmf}\\
{\cal R}({\bf x}_+):&={\cal S}_{k,n}^r\cap\{{\bf x}\,|\,{\bf 1}^{1\times r}{\bf x}={\bf x}_+\}\label{condpmf1};
\end{align}
$f_\text{mult}({\bf x}\,|\,{\bf x}_+)$ is a multivariate hypergeometric pmf. To show $\P_{\bf p}[L^*({\bf X})\ge c]$ is Schur-concave by applying the Proposition in Example 2 of [R] to the expectation in \eqref{PEP}, it must be shown that 
\begin{equation*}
\phi_c({\bf X}_+)\equiv\P_{\bf p}[L^*({\bf X})\ge c\,|\,{\bf X}_+]
\end{equation*}
 is Schur-concave in ${\bf X}_+$. Unfortunately this is not true in general: for the simplest case $k=n=r=2$, it can be shown that
\begin{equation}\label{222}
(\phi_c(0,4),\phi_c(1,3),\phi_c(2,2),\phi_c(3,1),\phi_c(4,0))=\ts(0,1,\frac{1}{3},1,0)
\end{equation}
when $\frac{1}{16}<c<\frac{2}{27}$, which violates Schur-concavity.

\vskip3pt

An alternative, more tractable procedure for testing problem \eqref{HvHm8} is obtained as follows. Factor the LRT statistic $L^*$ in \eqref{LRT9} as
\begin{align}
L^*({\bf X})&=\frac{\prod_{i=1}^r\prod_{j=1}^kX_{ij}!}{\prod_{j=1}^k X_{+j}!}\cdot\frac{\prod_{j=1}^k X_{+j}!}{\prod_{j=1}^kX_{+j}^{X_{+j}}}\nonumber\\
&\equiv V({\bf X}\,|\,{\bf X}_+)\cdot\~L({\bf X}+)\label{LRT11}\\
 &\propto\frac{f_\text{unif}({\bf X}\,|\,{\bf X}_+)}{f_\text{mult}({\bf X}\,|\,{\bf X}_+)}\cdot\~L({\bf X}+),\nonumber
\end{align} 
where we use the fact that the conditional pmf $f_\text{unif}({\bf X}\,|\,{\bf X}_+)$ remains uniform (constant) over the conditional range ${\cal R}({\bf x}_+)$. As noted above, $\~L({\bf X}_+)$ by itself may lose information relevant for \eqref{HvHm8}, but some of this missing information can be obtained from $V({\bf X}\,|\,{\bf X}_+)$ as now shown.

A simple but reasonable approach to combining $V({\bf X}\,|\,{\bf X}_+)$ with $\~L({\bf X}_+)$ for testing \eqref{HvHm8} is the following:
\vskip4pt
\nid($*$)\hskip60pt {\it reject $H_\mathrm{mult}$ if $V({\bf X}\,|\,{\bf X}_+)$ or $\~L({\bf X}_+)$ is large.} 
\vskip4pt

\nid As in \eqref{pLRT9}, we seek a tractable upper bound for the overall significance level. This can be accomplished by a hybrid procedure that combines the non-dependence on ${\bf p}$ of the conditional pmf $f_\text{mult}({\bf x}\,|\,{\bf x}_+)$ under  ${\cal M}(k,n)$ with the Schur-concavity in ${\bf p}$ of $\P_{\bf p}[\~L({\bf X}_+)\ge c]$. (Because ${\bf X}_+\sim\text{Mult}_{\bf p}(k,r n)$, Proposition 4.2 and Remark 4.3 apply to $\~L({\bf X}_+)$.) This procedure is now described.

For specified $\alp,\bet>0$,
 use $f_\text{mult}({\bf x}\,|\,{\bf x}_+)$ to determine $c_\alp({\bf x}_+)$ such that
\begin{equation}\label{calp}
\P[V({\bf X}\,|\,{\bf X}_+)\ge c_\alp({\bf X}_+)\,|\,{\bf X}_+]= \alp\ \ (\text{or}\approx\alp),
\end{equation}
and, as in \eqref{pLRT}, use $f_\mathrm{ecp}({\bf x})$ to determine $d_\bet$ such that
\begin{equation}\label{dalp}
\P_\text{ecp}[\~L({\bf X})\ge d_\bet]= \bet\ \ (\text{or}\approx\bet),
\end{equation}
Apply ($*$) to obtain the following test procedure for \eqref{HvHm8}:
\vskip4pt

\hskip40pt$T_{\alp,\bet}$: {\it reject $H_\mathrm{mult}$ if $V({\bf X}\,|\,{\bf X}_+)\ge c_\alp({\bf X}_+)$ or $\~L({\bf X}_+)\ge d_\bet$.} 
\vskip4pt

\nid{\bf Proposition 8.2.} The overall significance level of $T_{\alp,\bet}$ is $\le\alp+\bet$.
\vskip2pt

\nid{\bf Proof.}  By Bonferroni's inequality, the overall significance level of $T_{\alp,\bet}$ is
\begin{align}
&\sup_{{\bf p}\in{\cal P}_k}\P_{\bf p}[T_{\alp,\bet}\ \text{rejects}\ H_\mathrm{mult}]\nonumber\\
=&\sup_{{\bf p}\in{\cal P}_k}\P_{\bf p}[V({\bf X}\,|\,{\bf X}_{+})\ge c_\alp({\bf X}_+)\ \text{or}\ \~L({\bf X}_+)\ge d_\bet]\nonumber\\
\le&\sup_{{\bf p}\in{\cal P}_k}\{\P_{\bf p}[V({\bf X}\,|\,{\bf X}_{+})\ge c_\alp({\bf X}_+)]+\P_{\bf p}[\~L({\bf X}_+)\ge d_\bet]\}\nonumber\\
=&\;\alp+\sup_{{\bf p}\in{\cal P}_k}\P_{\bf p}[\~L({\bf X}_+)\ge d_\bet]\nonumber\\
 =&\;\alp+\P_\text{ecp}[\~L({\bf X}_+)\ge d_\bet]\nonumber\\
  =&\;\alp+\bet,\nonumber
\end{align}
which can be controlled by choosing $\alp$ and $\bet$ appropriately.\hfill$\square$


\newpage

\centerline{\bf Appendix 1}
\vskip10pt

\nid{\bf Proof of Proposition 5.2.} (i) The binomial pmf $f_p(x)$ is unimodal and
\begin{equation*}
\frac{n}{n+1}\le p\le1 \iff {n\choose n-1}p^{n-1}(1-p)\le{n\choose n}p,
\end{equation*}
so $f_p(x)$ is nondecreasing in $x$, hence the second sum in \eqref{ineq6.1} is
\begin{equation*}
\sum_{x=0}^n  (n-x+1)f_p(x)=n(1-p)+1.
\end{equation*}
Therefore the conjectured inequality \eqref{ineq6.1} is equivalent to
\begin{equation}\label{ineq11}
\sum_{x=0}^n |n-2x|f_{\frac{1}{2}}(x)> n(1-p).
\end{equation}
which, since $1-p\le\frac{1}{n+1}$, will hold if
\begin{equation}\label{ineq12}
\sum_{x=0}^n|n-2x|{n\choose x} > \frac{n2^n}{n+1}.
\end{equation}

Now apply the following identity, proved below:
\begin{equation}\label{lem1}
\sum_{x=0}^n |n-2x|{n\choose x}
=\begin{cases}n{n\choose\frac{n}{2}} &\text{if}\ n\ \text{is\ even},\\
 2n{n-1\choose\frac{n-1}{2}} &\text{if}\ n\ \text{is\ odd.}
\end{cases}
\end{equation}
If $n$ is even, it follows from this identity that \eqref{ineq12} will hold if
\begin{equation}\label{ineq13}
{n\choose\frac{n}{2}} > \frac{2^n}{n+1}.
\end{equation}
This inequality follows from a lower bound for ${n\choose\frac{n}{2}}$ in Krafft (2000):
\begin{equation*}
{n\choose\frac{n}{2}}\ge2^{n-1}\sqrt{\frac{2}{n}}> \frac{2^n}{n+1}.
\end{equation*}
Similarly, if $n$ is odd then \eqref{ineq12} will hold if
\begin{equation}\label{ineq14}
{n-1\choose\frac{n-1}{2}} > \frac{2^{n-1}}{n+1},
\end{equation}
which again follows from Krafft's inequality:
\begin{equation*}
{n-1\choose\frac{n-1}{2}}\ge2^{n-2}\sqrt{\frac{2}{n-1}}> \frac{2^{n-1}}{n+1}.
\end{equation*}

(ii) Recall that $q_{\frac{1}{2}}(x)=1+|n-2x|$ and $f_{\frac{1}{2}}(x)=f_{\frac{1}{2}}(n-x)$ for $x=0,\dots,n$. Furthermore,  if $p>\t12$ and $n/2<x\le n$ then
\begin{align}
\frac{f_p(x)}{f_p(n-x)}&=t^{2x-n}>1;\label{ratio}\\
\frac{f_p(x)}{f_p(n-x+1)}&=\frac{n-x+1}{x}t^{2x-n-1}
\begin{cases}<1&\text{if}\ t<a_n(x)\Leftrightarrow p<\frac{1}{2}+b_n(x),\\=1&\text{if}\ t=a_n(x)\Leftrightarrow p=\frac{1}{2}+b_n(x),\\>1&\text{if}\ t>a_n(x)\Leftrightarrow p>\frac{1}{2}+b_n(x),\end{cases}\label{ratio1}
\end{align}
where $t\equiv t(p)=\frac{p}{1-p}>1$ and $b_n(x)=\frac{a_n(x)-1}{2[a_n(x)+1]}$. 

First suppose that $n$ is odd. Then \eqref{ratio}, \eqref{ratio1},  and the unimodality of the binomial pmf imply that when $\t12< p<\t12+\eps_n$,
\begin{equation}\label{ratio2}
q_p(x)=\left.\begin{cases}q_{\frac{1}{2}}(x),&x=0,\dots,\frac{n-1}{2},\\q_{\frac{1}{2}}(x)-1,&x=\frac{n+1}{2},\dots,n\end{cases}\right\}.
\end{equation}
Thus in this case the conjectured inequality \eqref{ineq6.1} is equivalent to
\begin{align}
\sum_{x=0}^n |n-2x|f_{\frac{1}{2}}(x)&>\sum_{x=0}^n |n-2x|f_p(x)-\sum_{x=\frac{n+1}{2}}^nf_p(x)\label{ineq20}
\end{align}
Because
\begin{align}
\lim_{p\downarrow\frac{1}{2}}\left\{\sum_{x=0}^n |n-2x|f_p(x)-\sum_{x=\frac{n+1}{2}}^nf_p(x)\right\}=\sum_{x=0}^n |n-2x|f_{\frac{1}{2}}(x)-\t12,\label{ineq21}
\end{align}
the inequality \eqref{ineq20} holds when $\t12<p<\t12+\del_n$ for sufficiently small $\del_n<\eps_n$.

If $n$ is even then \eqref{ratio}, \eqref{ratio1},  and the unimodality of the binomial pmf imply that when $\t12< p<\t12+\eps_n$,
\begin{equation}\label{ratio2a}
q_p(x)=\left.\begin{cases}q_{\frac{1}{2}}(x),&x=0,\dots,\frac{n}{2},\\q_{\frac{1}{2}}(x)-1,&x=\frac{n}{2}+1,\dots,n\end{cases}\right\}.
\end{equation}
Thus in this case the conjectured inequality \eqref{ineq6.1} is equivalent to
\begin{align}
\sum_{x=0}^n |n-2x|f_{\frac{1}{2}}(x)&>\sum_{x=0}^n |n-2x|f_p(x)-\sum_{x=\frac{n}{2}+1}^nf_p(x)\label{ineq30}
\end{align}
Because
\begin{align}
\lim_{p\downarrow\frac{1}{2}}\left\{\sum_{x=0}^n |n-2x|f_p(x)-\sum_{x=\frac{n}{2}+1}^nf_p(x)\right\}=\sum_{x=0}^n |n-2x|f_{\frac{1}{2}}(x)-\t12\left[1-\frac{1}{2^n}{n\choose\frac{n}{2}}\right],\label{ineq31}
\end{align}
the inequality \eqref{ineq30} holds when $\t12p<\t12+\del_n$ for sufficiently small $\del_n<\eps_n$.\hfill$\square$
\vskip6pt

\nid{\bf Proof of identity \eqref{lem1}.} If $n\ge4$ is even then because ${n\choose x}={n\choose n-x}$, 
\begin{align*}
\sum_{x=0}^n |n-2x|{n\choose x}
&=\sum_{x=0}^{\frac{n}{2}-1} (n-2x){n\choose x}+\sum_{x=\frac{n}{2}+1}^{n} (2x-n){n\choose x}\\
&=2\left[ \sum_{x=\frac{n}{2}+1}^{n} x{n\choose x}-\sum_{x=0}^{\frac{n}{2}-1} x{n\choose x}\right]\\
&=2\left[ \sum_{x=0}^{\frac{n}{2}-1} (n-x){n\choose x}-\sum_{x=0}^{\frac{n}{2}-1} x{n\choose x}\right]\\
&=2\left[n \sum_{x=0}^{\frac{n}{2}-1} {n\choose x}-2\sum_{x=0}^{\frac{n}{2}-1} x{n\choose x}\right]\\
&=2\left[\frac{n}{2} \left\{\sum_{x=0}^{\frac{n}{2}-1} {n\choose x}+\sum_{x=\frac{n}{2}+1}^{n} {n\choose x}\right\}-2n\sum_{x=0}^{\frac{n}{2}-2} {n-1\choose x}\right]\\
&=2n\left[\frac{1}{2} \left\{2^n-{n\choose\frac{n}{2}}\right\}-\left\{2^{n-1}-{n-1\choose\frac{n}{2}-1}-{n-1\choose\frac{n}{2}}\right\}\right]
\end{align*}
\begin{align*}
&=2n\left[{n-1\choose\frac{n}{2}-1}+{n-1\choose\frac{n}{2}}-\frac{1}{2}{n\choose\frac{n}{2}}\right]\\
&=2n\left[{n\choose\frac{n}{2}}-\frac{1}{2}{n\choose\frac{n}{2}}\right]\\
&=n{n\choose\frac{n}{2}}.
\end{align*}

Next, if $n\ge3$ is odd then similarly to the even case,
\begin{align*}
\sum_{x=0}^n |n-2x|{n\choose x}
&=\sum_{x=0}^{\frac{n-1}{2}} (n-2x){n\choose x}+\sum_{x=\frac{n+1}{2}}^{n} (2x-n){n\choose x}\\
&=2\left[ \sum_{x=\frac{n+1}{2}}^{n} x{n\choose x}-\sum_{x=0}^{\frac{n-1}{2}} x{n\choose x}\right]\\
&=2\left[ \sum_{x=0}^{\frac{n-1}{2}} (n-x){n\choose x}-\sum_{x=0}^{\frac{n-1}{2}} x{n\choose x}\right]\\
&=2\left[n \sum_{x=0}^{\frac{n-1}{2}} {n\choose x}-2\sum_{x=0}^{\frac{n-1}{2}} x{n\choose x}\right]\\
&=2\left[\frac{n}{2} \left\{\sum_{x=0}^{\frac{n-1}{2}} {n\choose x}+\sum_{x=\frac{n+1}{2}}^{n} {n\choose x}\right\}-2n\sum_{x=0}^{\frac{n-3}{2}} {n-1\choose x}\right]\\
&=2n\left[\frac{1}{2} \left\{2^n\right\}-\left\{2^{n-1}-{n-1\choose\frac{n-1}{2}}\right\}\right]\\
&=2n{n-1\choose\frac{n-1}{2}}.
\end{align*}
This completes the proof of identity \eqref{lem1}.\hfill$\square$
\vskip8pt

\nid{\bf Proof of Lemma 6.1.} For $\frac{n}{2}+1\le x\le n-1$,
\begin{align*}
&(2x-n-1)[2(x+1)-n-1](D_{x+1}-D_{x})\\
=&\ (2x-n-1)\!\!\!\!\!\!\sum_{m=n-(x+1)+2}^{x+1} d_m-[2(x+1)-n-1]\!\!\!\sum_{m=n-x+2}^x d_m\\
=&\ (2x-n-1)\!\!\!\sum_{m=n-x+1}^{x+1} d_m-(2x-n+1)\!\!\!\!\!\!\sum_{m=n-x+2}^x d_m\\
=&\ (2x-n-1)(d_{n-x+1}+d_{x+1})-2\!\!\!\!\!\!\sum_{m=n-x+2}^x d_m\ \ >0,
\end{align*}
by the strict convexity of $d_m$ in $m$.\hfill$\square$
\vskip8pt

\nid{\bf Proof of Lemma 7.3.} (i) If no ties occur among $f_p(0),\dots,f_p(n)$ then $(r_p(0),\dots,r_p(n))$ is a permutation of $(0,\dots,n)$ so $\sum\limits_{x=0}^n  r_p(x)\~f_p(r_p(x))=\E(\~X_p)$.
Therefore (i) holds by \eqref{ftilderel}. When $\frac{n}{n+1}<p<1$, $f_p(n-1)<f_p(n)$ so by unimodality, no ties can occur.
\vskip2pt

(ii) Recall that $f_{\frac{1}{2}}(x)=f_{\frac{1}{2}}(n-x)$ and, by \eqref{qabs} and \eqref{rqrel},
\begin{align*}
r_{\frac{1}{2}}(x)&=n-|n-2x|=r_{\frac{1}{2}}(n-x).
\end{align*}
Thus, when $n$ is odd the tied pairs are $(0,n),(1,n-1),\dots,(\frac{n-1}{2},\frac{n+1}{2})$, so
\begin{align}
\sum_{x=0}^n r_{\frac{1}{2}}(x)f_{\frac{1}{2}}(x)&=2\sum_{x=0}^{\frac{n-1}{2}}r_{\frac{1}{2}}(x)f_{\frac{1}{2}}(x)=4\sum_{x=0}^{\frac{n-1}{2}}xf_{\frac{1}{2}}(x).\label{rfhalfsum}
\end{align}
Furthermore,
\begin{align*}
\~f_{\frac{1}{2}}(x)
&=\begin{cases}f_{\frac{1}{2}}(\frac{x}{2}),&x=0,2,4,\dots,n-1,\\f_{\frac{1}{2}}(\frac{x-1}{2}),&x=1,3,5,\dots,n,\end{cases}\end{align*}
so
\begin{align*}
\E(\~X_{\frac{1}{2}})&=\sum_{x=0}^nx\~f_{\frac{1}{2}}(x)\\
&\ts=\sum\limits_{x=0,2,\dots,n-1}xf_{\frac{1}{2}}(\frac{x}{2})+\sum\limits_{x=1,3,\dots,n}xf_{\frac{1}{2}}(\frac{x-1}{2})\\
&=\sum\limits_{x=0}^{\frac{n-1}{2}}2xf_{\frac{1}{2}}(x)+\sum\limits_{x=0}^{\frac{n-1}{2}}(2x+1)f_{\frac{1}{2}}(x)\\
&=4\sum_{x=0}^{\frac{n-1}{2}}xf_{\frac{1}{2}}(x)+\t12.
\end{align*}
By comparing this to \eqref{rfhalfsum}, we see that (ii) holds.

When $n$ is even, the tied pairs are $(0,n),(1,n-1),\dots,(\frac{n}{2}-1,\frac{n}{2}+1)$, so
\begin{align}
\sum_{x=0}^n r_{\frac{1}{2}}(x)f_{\frac{1}{2}}(x)&=2\sum_{x=0}^{\frac{n}{2}-1}r_{\frac{1}{2}}(x)f_{\frac{1}{2}}(x)+\ts nf_{\frac{1}{2}}(\frac{n}{2})=4\sum\limits_{x=0}^{\frac{n}{2}-1}xf_{\frac{1}{2}}(x)+\ts nf_{\frac{1}{2}}(\frac{n}{2}).\label{rfhalfsum1}
\end{align}
Furthermore,
\begin{align*}
\~f_{\frac{1}{2}}(x)
&=\begin{cases}f_{\frac{1}{2}}(\frac{x}{2}),&x=0,2,4,\dots,n-2,n,\\f_{\frac{1}{2}}(\frac{x-1}{2}),&x=1,3,5,\dots,n-1,\end{cases}\end{align*}
so
\begin{align*}
\E(\~X_{\frac{1}{2}})&=\sum_{x=0}^nx\~f_{\frac{1}{2}}(x)\\
&\ts=\sum\limits_{x=0,2,\dots,n-2}xf_{\frac{1}{2}}(\frac{x}{2})+\sum\limits_{x=1,3,\dots,n-1}xf_{\frac{1}{2}}(\frac{x-1}{2})+ \ts nf_{\frac{1}{2}}(\frac{n}{2})\\
&=\sum\limits_{x=0}^{\frac{n}{2}-1}2xf_{\frac{1}{2}}(x)+\sum\limits_{x=0}^{\frac{n}{2}-1}(2x+1)f_{\frac{1}{2}}(x)+ \ts nf_{\frac{1}{2}}(\frac{n}{2})\\
&=4\sum_{x=0}^{\frac{n}{2}-1}xf_{\frac{1}{2}}(x)+\t12+\ts nf_{\frac{1}{2}}(\frac{n}{2}).
\end{align*}
By comparing this to \eqref{rfhalfsum1}, we see that (ii) again holds.\vskip2pt

(iii) If ties occur among $f_p(0),\dots,f_p(n)$ then by the unimodality of $f_p(x)$ in $x$, these ties must occur in nested non-overlapping pairs $f_p(y_i)=f_p(x_i)$, where $1\le\cdots<y_2<y_1<x_1<x_2<\cdots\le n$. For fixed $n$, the most such tied pairs occur when $p=\t12$ as in (ii). Here we consider the case $\t12<p\le\frac{n}{n+1}$.


Partition the half-open interval $(\t12,\frac{n}{n+1}]$ as follows:
\begin{align*}
\left(\t12,\frac{n}{n+1}\right]=
\begin{cases}\ts\bigcup_{x=\frac{n+1}{2}}^{n-1}\left(\frac{x}{n+1},\frac{x+1}{n+1}\right]&\text{if}\ n\ \text{is\ odd},\\
\ts\left(\t12,\frac{\frac{n}{2}+1}{n+1}\right]\cup\bigcup_{x=\frac{n}{2}+1}^{n-1}\left(\frac{x}{n+1},\frac{x+1}{n+1}\right]&\text{if}\ n\ \text{is\ even}.
\end{cases}
\end{align*}
Because
\begin{equation*}
\ts p\in\left(\frac{x}{n+1},\frac{x+1}{n+1}\right]\implies f_p(x-1)<f_p(x)\ge f_p(x+1),
\end{equation*}
by unimodality, ties can occur only for $(y_i,x_i)$ pairs with $x+1\le x_i\le n$.  For each such tie, the coefficient of $f_p(x_i)$ in 
$\E(\~X_p)$ is greater by 1 than its coefficient $r_p(x_i)$ in $\sum_{x=0}^n  r_p(x)f_p(x)$. Therefore, since $p\le\frac{x+1}{n+1}$ and the binomial$(n,p)$ distribution is stochastically increasing in $p$, 
\begin{align}
\Del_{n,p}:&=\E(\~X_p)-\sum_{z=0}^n  r_p(z)f_p(z)
\le\sum_{z=x+1}^nf_p(z)\le\sum_{z=x+1}^nf_{\frac{x+1}{n+1}}(z).\label{Del}
\end{align}
However, by the relationship between the binomial and beta distributions,
\begin{align}
\sum_{z=x+1}^nf_{\frac{x+1}{n+1}}(z)&=\ts\Pr[\text{Beta}(x+1,n-x)\le\frac{x+1}{n+1}]<\t12.\label{Del1}
 \end{align}
This inequality holds because median[Beta$(\alp,\bet)]>\frac{\alp}{\alp+\bet}$ when $\alp>\bet>1$ (cf. Kerman (2011)), hence $\Del_{n,p}<\t12$. Since $\Del_{n,p}>0$ because at least one tie is present, (iii) is verified for $p\in\left(\frac{x}{n+1},\frac{x+1}{n+1}\right]$.
 
Lastly, because
\begin{equation*}
\ts p\in\left(\t12,\frac{\frac{n}{2}+1}{n+1}\right]\implies f_p(\frac{n}{2}-1)<f_p(\frac{n}{2})\ge f_p(\frac{n}{2}+1),
\end{equation*}
by unimodality, ties can occur only for $(y_i,x_i)$ pairs with $\frac{n}{2}+1\le x_i\le n$.  
Since $p\le\frac{\frac{n}{2}+1}{n+1}$, it follows as in \eqref{Del} and \eqref{Del1} with $x=\frac{n}{2}$ that
\begin{align*}
\Del_{n,p}&\le\sum_{z=\frac{n}{2}+1}^nf_{\frac{\frac{n}{2}+1}{n+1}}(z)<\t12.
\end{align*}
Again $\Del_{n,p}>0$ because at least one tie is present, so (iii) is verified for $p\in\left(\t12,\frac{\frac{n}{2}+1}{n+1}\right]$.
 \vskip2pt

(iv) If $p=1$ then $n-1$ ties occur, namely
\begin{align*}
f_1(0)=\cdots=f_1(n-1)=0,&\quad f_1(n)=1,\\
\~f_1(0)=\cdots=\~f_1(n-1)=0,&\quad \~f_1(n)=1,\\
r_1(0)=\cdots=r_1(n-1)=0,&\quad r_1(n)=n,
\end{align*}
so (iv) holds, since both sides of (iv) equal $n$.\hfill$\square$
\vskip12pt

\vskip14pt


\centerline{\bf Appendix 2}
\vskip10pt

\nid{\bf A critique of pure significance tests.} The PST method described in Section 1 possesses a seemingly irreparable ambiguity, as now demonstrated.
\vskip4pt

\nid{\bf Example 1.} Suppose that Mike observes the rv $U$ whose range is the interval $(0,1)$, Joe observes $X=-\log U$ whose range is $(0,\infty)$, and Steve observes $Y=-\log(1-U)$ whose range also is $(0,\infty)$. Note that $U$, $X$, and $Y$ are mutually equivalent, so Mike, Joe, and Steve possess the same information. 

Suppose that Mike wishes to test $H_0:U\sim\mathrm{Uniform}(0,1)$, equivalently, $f_U(u)=1_{(0,1)}(u)$, the uniform pdf on $(0,1)$, without specifying any alternative distribution(s). Because $f_U$ is constant, every PST based on $U$ is trivial, either rejecting $H_0$ for all $U$ or accepting $H_0$ for all $U$. 

However, the situation is much worse for Joe and Steve. It is straightforward to show that for Joe, $H_0$ is equivalent to $f_X(x)=e^{-x}1_{(0,\infty)}(x)$, the standard exponential pdf, while for Steve, $H_0$ is equivalent to $f_Y(y)=e^{-y}1_{(0,\infty)}(y)$, also standard exponential. Thus for Joe the PST rejects $H_0$ if $e^{-X}\le\tau$ where $0<\tau<1$, equivalently, if
\begin{equation}\label{Joe}
X\ge-\log\tau.
\end{equation}
Similarly, for Steve the PST rejects $H_0$ if $e^{-Y}\le\tau$, equivalently, if $Y\ge-\log\tau$. However,
\begin{equation*}
Y=-\log(1-U)=-\log(1-e^{-X}),
\end{equation*} 
so Steve's PST rejects $H_0$ if 
\begin{equation}\label{Steve}
X\le-\log(1-\tau),
\end{equation}
which is essentially the opposite of Joe's PST \eqref{Joe}. In fact, in terms of $U$, \eqref{Joe} and \eqref{Steve} become, respectively,
\begin{align}
U&\le\tau,\label{Joe1}\\
U&\ge1-\tau.\label{Steve1}
\end{align}
Thus if $\tau=\t12$ then \eqref{Joe} and \eqref{Steve} become
\begin{align}
U&\le\t12,\label{Joe2}\\
U&\ge\t12,\label{Steve2}
\end{align}
so the two PSTs reach {\it exactly} opposite conclusions, although based on equivalent evidence. Clearly this is an undesirable property for an inference procedure.

And the situation is actually far worse. Suppose that Marina observes $Z=-\log\eta(U)$, where $\eta$ is an {\it arbitrary} measure-preserving bijection of $(0,1)\to(0,1)$. Under $H_0$, $\eta(U)\sim U$, so for Marina $H_0$ is equivalent to $f_Z(z)=e^{-z}1_{(0,\infty)}(z)$, again standard exponential. Thus her PST rejects $H_0$ if $Z\ge-\log\tau$, equivalently, if
\begin{equation}\label{Marina}
\eta(U)\le\tau.
\end{equation}
Because this holds for all $\eta$, this shows that in terms of $U$ {\it every} measurable subset of $(0,1)$ is a rejection region for a PST of $H_0$. Thus, in terms of $Z$ essentially {\it every} measurable subset of $(0,\infty)$ is the rejection region of a PST for $H_0$, again an unacceptable conclusion.\hfill$\square$
\vskip4pt

\nid{\bf Example 2.} For the {\it coup de grace}, we show that the above discussion extends to essentially {\it all} pdfs $f_0$ on $\mathbb{R}^1$, not just standard exponentials. Suppose it is wished to test $H_0:f_X=f_0$ without specifying any alternative(s). For simplicity, assume that the support of $f_0$ is an interval $(a,b)$, where $-\infty\le a<b\le\infty$ and let $F_0$ be the cdf corresponding to $f_0$. Under $H_0$, $U:=F(X)\sim\mathrm{Uniform}(0,1)$. ($U$ is the probability integral transform of $X$.)
The rejection region of the PST for $H_0$ is
\begin{align}
R_\tau&=\{x\in(a,b)\,|\,f_0(x)\le\tau\}\label{R}\\
 &=\{u\in(0,1)\,|\,f(F_0^{-1}(u))\le\tau\}.\nonumber
\end{align}

Now define
\begin{equation}\label{Reta}
W_\eta=F_0^{-1}(\eta(F_0(X)))=F_0^{-1}(\eta(U)),
\end{equation}
where again $\eta$ is an {\it arbitrary} measure-preserving bijection of $(0,1)\to(0,1)$. It is straightforward to show that $f_{W_\tau}=f_0$ under $H_0$, so the PST for $H_0$ based on $W_\tau$ has rejection region
\begin{align}
\{w\in(a,b)\,|\,f_0(w)\le\tau\}&=\{u\in(0,1)\,|\,f_0(F_0^{-1}(\eta(u)))\le\tau\}\label{Rtau1}\\
 &=\eta^{-1}\left(\{v\in(0,1)\,|\,f_0(F_0^{-1}(v))\le\tau\}\right)\nonumber\\
  &=\eta^{-1}(R_\tau).\nonumber
\end{align}
Because this holds for all $\eta$, in terms of $U$ essentially {\it every} measurable subset of $(0,1)$ is the rejection region of a PST for $H_0$. Thus, in terms of $X$ essentially {\it every} measurable subset of $(a,b)$ is the rejection region of a PST for $H_0$, again an undesirable property.\hfill$\square$

\vskip15pt

\centerline{\bf References}

\def\new{\hangindent=1.5pc}

\vskip .25cm \new
\noindent [C] Chernoff, H. (1956). Large sample theory: parametric case. {\it Ann. Math. Statist.}  {\bf 27} 1-22.

\vskip .25cm \new
\noindent [H] Hodges, J. S. (1990). Can/may Bayesians do pure tests of significance? In {\it Bayesian and Likelihood Methods in Statistics and Econometrics}, S. Geisser, J. S. Hodges, S. J. Press, A. Zellner, eds., 75-90. Elsevier (North Holland).

\vskip .25cm \new
\noindent [H] Howard, J. V. (2009). Significance testing with no alternative hypothesis: a measure of surprise. {\it Erkenntnis} {\bf 70} 253-270.

\vskip .25cm \new
\noindent [K] Kerman, J. (2011). A closed-form expression for the median of the beta distribution. arXiv:1111.0433.

\vskip .25cm \new
\noindent [K] Krafft, O. (2000).  Problem 10819, {\it Amer. Math. Monthly} {\bf 107} p.652.

\vskip .25cm \new
\noindent [MOA] Marshall, A. W, I. Olkin, B. C. Arnold (2011).  {\it Inequalities: Theory of Majorization and its Applications (2nd ed.)}, Springer, New York.

\vskip .25cm \new
\noindent [R] Rinott, Y. (1973). Multivariate majorization and rearrangement inequalities with some applications to probability and statistics. {\it Israel J. Math.}  {\bf 15} 60-77.

\vskip .25cm \new
\noindent [SS] Sackrowitz, H. and E. Samuel-Cahn (1999).  $P$ values as random variables - expected $p$ values. {\it Amer. Statistician} {\bf 53} 326-331.

\end{document}